\documentclass[12pt]{article}
 \usepackage{amsfonts,amssymb,graphicx}
\usepackage{inputenc}
\usepackage[pdfpagemode=None,colorlinks=true,linkcolor=blue,citecolor=blue]{hyperref}

\newtheorem{prop}{Proposition}
\newtheorem{rema}{Remark}

\newtheorem{lemm}{Lemma}
\newtheorem{theo}{Theorem}
\newtheorem{coro}{Corollary}

\newcommand{\R}[1][]{\ensuremath{{\mathbb{R}^{#1}} }}
\renewcommand{\S}[1][]{\ensuremath{{\mathbb{S}^{#1}} }}

\newcommand{\<}{\langle}
\renewcommand{\>}{\rangle}

\newcommand{\pa}{\partial}

\newcommand{\al}{\alpha}
\newcommand{\eps}{\epsilon}

\newcommand{\be}{\beta}

\newcommand{\w}{\tilde{w}}
\date{}

\title{On the three-dimensional Blaschke-Lebesgue problem}
\author{Henri Anciaux\footnote{The first author is supported by Science Foundation Ireland (Research Frontiers Program)} , Brendan Guilfoyle}
\begin{document}

\maketitle

\centerline{\textbf {Abstract}}

\smallskip

{\small The width of a closed convex subset of $n$-dimensional Euclidean space is the distance between two
parallel supporting hyperplanes. The  Blaschke-Lebesgue problem consists
 of minimizing the volume in the class of convex sets of fixed constant width and is still open in dimension $n \geq 3.$ 
In this paper we  describe a necessary condition that  the minimizer of the 
Blaschke-Lebesgue must satisfy in dimension $n=3$: we prove that the smooth components of the boundary of the minimizer
have their smaller principal curvature constant, and therefore
are either spherical caps or pieces of tubes (canal surfaces).}

\medskip

\centerline{\small \em 2000 MSC: 52A40,  52A15
\em }


\section*{Introduction}
The width of a convex body $B$ in $n$-dimensional Euclidean
space in the direction $\vec{u}$ is the distance between the two
supporting hyperplanes of $B$ which are orthogonal to $\vec{u}.$ When
this distance is independent of $\vec{u}$, $B$ is said to have \em
constant width. \em The ratio ${\cal I}(B)$ of the volume of a constant width body to 
the volume of the ball of the same width
 is a homothety invariant, as is the isoperimetric ratio. Moreover, the maximum of
 ${\cal I }(B)$, like  the minimum of the isoperimetric ratio, is
 attained by  round spheres.
 However, while the isoperimetric ratio is not bounded from above, the infimum of ${\cal I}$ is strictly
 positive, since
compactness properties of the space of convex sets ensures the
existence of a minimizer.  It is known
by the work of Blaschke and Lebesgue that the \em Reuleaux triangle, \em
obtained by taking the intersection of three discs centered at the
vertices of an equilateral triangle, minimizes ${\cal I}$ in
dimension $n=2.$ The determination of the minimizer of ${\cal I}$ in any dimension is
referred to as
the  \em Blaschke-Lebesgue problem. \em 

\medskip

Recently several simpler solutions of the problem in dimension $2$
have been given (see \cite{Ba},\cite{Ha}), however the Blaschke-Lebesgue
problem in dimension $n \geq 3$ appears to be very difficult to solve
and remains open. A crucial step in solving  the Blaschke-Lebesgue problem in dimension $n=2$ consists of 
proving that the boundary of the minimizer is made up of arcs of circles of radii equal to the width,
 and hence the smooth parts of the boundary 
have constant curvature.

\medskip

In this paper we give a property of the minimizer of the Blaschke-Lebesgue in dimension $n=3$  which generalizes
the constant curvature condition observed in dimension $n=2$ (here and in the following, "smooth" means "continuously twice differentiable"):

\bigskip

\noindent \textbf{Main Theorem:} \em Let  $B$ be a local  minimizer of the Blaschke-Lebesgue problem in $\mathbb R^3$
 with constant width $2w.$
 
Then the smooth parts of its boundary  have their smaller principal curvature  constant and equal to $1/2w.$ \em

\bigskip

It is easily seen that the boundary of a constant width body in $\R^3$ cannot be made up only of spherical caps, so 
the minimizer of the Blaschke-Lebesgue problem must have a more complicated geometry. On the other hand, 
K. Shiohama and R. Takagi 
 proved in \cite{ST} that a non-spherical surface with one constant principal curvature must be a canal surface, i.e.\
 the envelope of a one-parameter family of spheres, or, equivalently, a tube over a curve (i.e.\ 
  the set of points which lie at a fixed distance from this curve).
  Thus the main theorem implies the following:

\bigskip

\noindent \textbf{Corollary:}
\em Let $B$ be a local  minimizer of the Blaschke-Lebesgue problem in $\mathbb R^3$ with constant width $2w.$

Then the smooth parts of its boundary are spherical caps or pieces of tubes, both of them with radius equal to the width $2w$ of $B.$ \em

\bigskip

We observe that the constant width body having the best known ratio ${\cal I}$,
 Meissner body (\cite{CG}, \cite{GK}, \cite{Ba}) satisfies this criteria: it is made up of 
 four spherical caps centered at the vertices of a tetrahedron, and three tubes over three arcs of
 circles. Therefore
we cannot discard the possibility that it is the solution of the Blaschke-Lebesgue problem, although 
one might expect the minimizer to  have tetrahedral symmetry. Another interesting
constant width body is the one obtained by rotation of the Reuleaux triangle about an axis of symmetry.
It is known that the latter minimizes the ratio ${\cal I}$ among constant width bodies with rotational symmetry
(see \cite{CCG}, \cite{AG}). It is also interesting to note that this body  satisfies our criteria as well:
 one part of its boundary is a spherical cap, and the other one is a tube
over an arc of a circle. However, it has a bigger ratio ${\cal I}$ than Meissner's, which in particular
proves that the solution of the Blaschke-Lebesgue problem does not have rotational symmetry.

\medskip

In the light of our result,  the most difficult issue to address seems to be that of the regularity.  
We cannot exclude
\em a priori \em that the boundary of the minimizer of the Blaschke-Lebesgue problem is singular everywhere and
the traditional techniques of regularity theory (e.g.\ those used for harmonic maps or minimal surfaces) do not seem to apply here.
On the other hand, assuming that
the minimizer is made up of a finite number of smooth parts, our result reduces the problem  
to a kind of combinatoric (though not easy) one: minimize the volume among the convex bodies whose boundaries are 
made up of spherical caps  and pieces of tubes, all of them of the same radius.

\medskip

As in \cite{Ba}, \cite{Ha} and \cite{AG}, our proof is based on the analysis
of the support function $s$ which characterizes a convex body $B$ of
constant width $2w.$ The first point consists of evaluating the volume of $B$ and the area
of its boundary in terms of $w$ and the function $h=s-w$ (Theorem \ref{VA}). Our formula allows us
in particular to prove easily the famous Blaschke formula, a functional relation between
the volume, the area and the width of $B,$ and to recover 
the fact that the ratio ${\cal I}$ 
is maximized by  round spheres.
A crucial point is then the following
observation, stated in \cite{GK}: flowing the boundary of a convex
body along its inward unit normal vector field preserves the
constant width condition, as long as the evolving surface remains
convex.
 Moreover, the ratio ${\cal I}$ decreases along the flow, so
the minimizer of ${\cal I}$ must occur at the latest time such
that convexity holds, and therefore must be singular. This issue
is easily controlled since the function $h$ 
 is invariant along the normal flow, while the width
$2w$ decreases linearly. Thus, there
exists a positive number $w_0(h)$ such that for any $w \geq
w_0(h),$  the function $s=h+w$ is the support function  of some convex body of
constant width $2w$.
 Hence, we can restrict the minimization process to the class of support functions of the form $s=h+w_0(h),$
while all the necessary information is carried by the function
$h.$ The main theorem is then obtained as follows: assuming that the smaller principal curvature is not
constant on some smooth part of the boundary, we compute the second variation of ${\cal I}$ 
for a suitable local deformation of $h$ to get a contradiction.

\medskip

The authors wish to thank the referee for correcting the statement of Theorem~\ref{const}.

\section{The geometry of constant width bodies}
Let $B$ be a convex body in $\R^{n}$ and denote by $s$ its support function, i.e.\
$ s(u)=\sup_{x \in B} \<u,x \> , \forall u \in \S^{n-1}.$
Then the width $w(u)$ of $B$ in the direction $u$ is related to the support function by the following formula:
 $$ 2w(u)= s(u)+ s( -u),$$
 where $-u$ is the antipodal point of $u$ in $\S^{n-1}.$
It is known (see \cite{Ho}, \cite{Ba}) that if $B$ has constant width it must be
strictly convex; moreover  it is proven in \cite{Ho} that the support function $s$ of a constant width body is $C^{1,1},$
i.e.\ it admits first derivatives which are Lipschtiz continuous. By the Rademacher theorem, it follows that the
second derivatives are well defined almost everywhere and bounded. This fact will be important later on; since the geometry
of the boundary of $B$ will be expressed in terms of the Hessian of $h.$
  
If $B$ is a strictly convex body in $\R^n$ whose support function $s$ belongs to $C^{1,1},$ the following map
$$ \begin{array}{lccc} f :
 &  \mathbb S^{n-1}  &\to& \mathbb R^{n}              \\
&  u & \mapsto & s(u).u+ \nabla s(u)
   \end{array}$$
 is a parametrization of its boundary and $u$ is the Gauss map of $\partial B$.

Given an arbitrary strictly convex body $B$, let $w \in \R$ 
be the mean of its support function $s$ on $\S^{n-1}$:
$$ w:=\frac{\int_{\mathbb S^{n-1}} s(u) dA}{\int_{\mathbb S^{n-1}} dA},$$
where $dA$ denotes the canonical volume form on $\S^{n-1},$ and introduce the zero mean map
$ h:=s-w.$ 

Then $B$ has constant width $2w$ if and only if the function $h$ is odd, i.e.\
$$h(u)+h(-u)=0.$$ 

The following inequality will be crucial for us:
\begin{prop}[Wirtinger inequality] \label{Wirtinger}
Let $h \in C^{1,1}(\mathbb S^{n-1})$ with vanishing mean and $dA$ the volume element on $\mathbb S^{n-1}$. Then
$${\cal E}(h):=\int_{\mathbb S^{n-1}} \left(\frac{1}{n-1}|\nabla h|^2-h^2\right) dA \geq 0,$$
 with equality if $h$ is a first eigenfunction of the Laplacian
on the sphere $\mathbb S^{n-1}$.
\end{prop}

\noindent This result is easily proved once the theory of spherical harmonics, generalizing
  Fourier analysis to higher dimension, is developed (see  \cite{GW}, p.\ 1288).

\bigskip

From now on, we restrict ourselves to the case of dimension $3.$ Our first step consists of expressing
the local geometry of the boundary of a convex body $B$ in terms of the data $(h,w)$. We recall that the Hessian
of $h$ is the symmetric tensor defined by $Hess(h) (X,Y)=\<\nabla_X \nabla h,Y\>,$ where $\nabla$ denotes the
Levi-Civita connection of the round metric of $\S^2.$ The two invariants of $Hess(h)$ are
its trace, which is the well known Laplace-Beltrami operator $\Delta$ and its determinant, that we shall denote
in the following by $H(h).$

\begin{theo} \label{geo}
The area element of $\pa B$, denoted by $d\bar{A}$, is given by:
$$ d\bar{A}=(w^2 + \al w + \be)dA,$$
where we set
$$\al:=2h+\Delta h  \quad \mbox{ and } \quad \be:=h^2+h\Delta h+ H(h) .$$
Moreover, its principal curvatures $k_1$ and $k_2$, whenever they exist, take the following form:
$$ k_{1,2} = \frac{2w + \al \pm \sqrt{\al^2-4 \be}}{2(w^2 + \al w+ \be)}.$$
\end{theo}

  In the case where $B$ has constant width, we deduce the following formulas for its volume
 ${\cal V}(B)$  and the area of its boundary ${\cal A}(\pa B)$:

\begin{theo} \label{VA} Let $B$ be a convex body of constant width $2w$ in $\Bbb R^3.$ Then:
$${\cal V}(B)= \frac{4 \pi}{3} w^3 -w{\cal E}(h),$$
$${\cal A}(\pa B)= 4 \pi w^2  - {\cal E}(h).$$
 \end{theo}

This allows us to recover the famous Blaschke formula, a functional relation between the volume, the area and the width:

\begin{coro}[Blaschke formula] Let $B$ be a convex body in $\Bbb R^3$ of constant width $2w.$ Then:
$${\cal V}(B)=w {\cal A}(\partial B) - \frac{8}{3} \pi w^3.$$
\end{coro}

The proofs of Theorem \ref{geo} and \ref{VA} are postponed at the end of the paper (sections 4 and 5).

\section{The Blaschke-Lebesgue problem}
Let $B$ be a convex body of constant width $2w$ and denote by $B_w$ the round ball  of radius $w.$
Introduce the ratio
$${\cal I}(B)= \frac{{\cal V}(B)}{{\cal V}(B_w)}=\frac{{\cal V}(B)}{4 \pi w^3/3}.$$
By Theorem \ref{VA}, we have
$${\cal I}(B)={\cal I}(h,w)= 1 - \frac{{\cal E}(h)}{4 \pi w^2/3}.$$
It follows from the Wirtinger inequality that 
the ratio ${\cal I}(B)$ is less than or equal to  $1$ and  equality is attained  when
 $h$ is a first eigenfunction of the Laplacian, as it is the case of balls $B=B_w$. 
 Moreover, for a given $h$,  ${\cal I}$ increases with respect to
$w$.  Hence it reaches its minimum at the lowest value of $w$ such that $h+w=s$
is the support function of a convex body; we define $w_0(h)$ to be this crucial quantity. 
Increasing (resp.\ decreasing) the value of $w$ corresponds geometrically to flowing
 the boundary of $B$ parallel to itself, i.e.\
 along its outward (resp.\ inward) normal vector. Therefore the map $h$ corresponds to a one-parameter family of parallel surfaces, labelled
 by the parameter $w \in [w_0(h),\infty).$
 The inward  normal flow can be continued
as long as the surface is smooth. By Theorem \ref{geo}, this is equivalent to the fact that
the area element $d\bar{A}$ is strictly positive.
 Hence,  we deduce
an explicit expression for $w_0(h)$:
$$w_0(h)=\inf \left\{ w \in \R^+ | \, \,  w^2 + \al w + \be > 0 \mbox{ a.e. on } \S^2  \right\}$$ 
and the convex body $B$ corresponding to $s=h+w_0(h)$ is always singular.
\begin{rema} 
One can check that $w_0(h)=||  W(h)||_{L^{\infty}(\S^2)}, $
where
$$ W(h)(u):=   \frac{-\al + \sqrt{\al^2-\be}}{2}.$$ 
\end{rema}

The directions $u$ of $\S^2$ where the area element vanishes corresponds precisely to points $f(u)$ of the boundary which are singular.
 The next theorem shows that in the smooth parts of a local minimizer of
${\cal I},$  such a situation actually occurs for  every pair of antipodal directions $(u,-u)$. We point out that this result is roughly equivalent to one of the main results of \cite{BLO} (Theorem 5).

\begin{theo} \label{const}
Let $(h,w_0(h))$ be a local minimizer of ${\cal I}(h,w)$ and let $U$ be an open subset of $\S^2$ where $h$ is smooth.
Then for every point $u$ of $U,$ the
area element $d\bar{A}$ vanishes 
 at one of the points
$u$ and $-u$.
\end{theo}

\noindent \textit{Proof.} We proceed by contradiction assuming
that there is an open subset $U$ of $\S^2$ where $h$ is smooth and 
such that $(w_0(h))^2 + \al w_0(h) + \be > 0$  in $U \cup (-U)$.
  Consider a smooth map $v$ such that
  $v(u)+v(-u)=0, \forall u \in \S^2$ and
   whose support
is contained in $U \cup (-U)$  
and define the deformation $h^\eps:=h+\eps v $ of $h.$ For
small $\eps,$
$$w_0(h^\eps)=w_0(h),$$
hence
$$\frac{{\cal E}(h^\eps)}{w^2_0(h^\eps)}  =\frac{{\cal E}(h)}{w^2_0(h)} + \eps  \frac{\delta {\cal E}(h,v)}{w^2_0(h)}
+ \frac{\eps^2}{2} \frac{ \delta^2 {\cal E}(h,v)}{w^2_0(h)} +
o(\eps^2).$$ As $h$ is a minimizer of ${\cal I}$, and thus a maximizer of
${\cal E}(h)/w^2_0(h)$, we must have both $\delta {\cal E}(h,v)
=0$ and $\delta^2 {\cal E}(h,v) \leq 0.$ On the other hand the
functional ${\cal E}$ is quadratic, so that $\delta^2 {\cal
E}(h,v)= {\cal E}(v)$, which is positive by the Wirtinger inequality
(Proposition \ref{Wirtinger}). Finally, the support of $v$ being contained in $U\cup (-U)$, $v$ cannot be an 
eigenfunction of the Laplacian and we get the required contradiction.

\section{Proof of the main theorem} \label{geometric}

We are now in position to prove our main result: assume that $B$ is a local minimizer of ${\cal I}(B)$
and let $h$ be the associated map. For the sake of brevity we set $\w:=w_0(h)$ in the following.
Let $U$ be an open subset 
of $\S^2$ such that $f(U)$ is a smooth part of $\pa B.$ In particular, by Theorem \ref{geo}, 
$\w^2 + \al \w + \be >0$ on $U.$
Hence by Theorem \ref{const}, 
$  \w^2 + \al(u)\w + \be(u)=0, \forall u \in -U.$
Since $\al$ is odd and $\be$ is even, it follows that
$  \w^2 - \al(u)\w + \be(u)=0,$ $\forall u \in U.$ 
Consequently, by Theorem \ref{geo},
$$ k_{1,2} =\frac{2\w + \al \pm \sqrt{\al^2-4 \be}}{2(\w^2 + \al \w+ \be)}=
\frac{2\w + \al \pm \sqrt{\al^2-4\al \w +4 \w ^2}}{2(\w^2 + \al \w+ (\al \w - \w^2))}$$
$$=\frac{2\w + \al \pm |\al - 2 \w|}{4 \al \w},$$
so that $k_1 = \frac{1}{\al}$ and $k_2=\frac{1}{2\w}.$
Hence the principal curvature $k_2$ is constant  on $U$ and equal to the inverse of the width $2\w$ of $B.$
 Finally, since $\al(u) \geq 0,$ we have
$$2\w={-\al(-u) + \sqrt{\al^2(-u) -4\be(-u)}}=\al(u) + \sqrt{\al^2(u) -4\be(u)} \geq {\al},$$
 so $ \frac{1}{2\w} \leq \frac{1}{\al}$, that is $k_2$ is the smaller curvature and the proof is complete.

\section{The local geometry of $\pa B$ (proof of Theorem \ref{geo})}
Let $(x,y)(u)$ be an isothermic coordinate chart from a dense subset  of 
$\S^2$ into an open subset $U$ of $\R^2$ (for example the
stereographic projection onto $U=\R^2$) and
denote by $e^r$ the conformal factor, i.e.\ $e^r=|\pa_x|=|\pa_y|.$ In particular
the area element is given by $dA= e^r dx dy.$ 
The coefficients of the Hessian of $h$ in the coordinates $(x,y)$ are:
$$a:=e^{-2r}\<\nabla_{\pa_x} \nabla h, \pa_x\>,$$
$$b:=e^{-2r}\<\nabla_{\pa_x} \nabla h, \pa_y\>=e^{-2r}\<\nabla_{\pa_y} \nabla h, \pa_x\>,$$
$$c:=e^{-2r}\<\nabla_{\pa_y} \nabla h, \pa_y\>.$$
We recall that the boundary of $\pa B$ is parametrized by $f(u)=s(u)u+ \nabla s(u)=s(u)u+ \nabla h(u).$
In order to compute the first derivatives of $f,$ we use the Gauss formula of the embedding of the sphere $\S^2$
in $\R^3,$ which relates the flat connection $D$ of $\R^3$ to the Levi-Civita connection $\nabla$ on the sphere:
$$(D_X Y)(u) = (\nabla_X Y)(u) - \<X,Y\>u.$$
 It follows that
 $$ f_x=s \pa_x + \nabla_{\pa_x} \nabla h=(s+a)\pa_x+b\pa_y,$$
and
 $$ f_y =s \pa_y + \nabla_{\pa_y} \nabla h=b\pa_x+(s+c)\pa_y.$$
The 
 trace $a+c=\Delta h$ and the determinant $H(h):=ac-b^2$ of the Hessian matrix of $h$
 are 
intrinsic quantities, i.e.\ they depend only on the metric on $\S^2,$ and not on the choice of coordinates.
 
We then compute the coefficients of the first fundamental form of the immersion~$f$:
$$E:= \<f_x,f_x\> = ((s+a)^2 +b^2)e^{2r},  \quad \quad F:=\<f_x,f_y\>= (2s+a+c)be^{2r}, $$
$$G:=\<f_y,f_y\>=((s+c)^2 +b^2)e^{2r}.$$
It follows that
$$\sqrt{EG-F^2}e^{-2r}=\Big( ((s+a)^2 +b^2)((s+c)^2 +b^2)-4(s+a+s+c)^2 b^2\Big)^{1/2} $$
$$=\left( s^2 + s(a+c)+ac-b^2\right)=\left( w^2 + (2h +a+c)w+ h^2 + (a+c)h +ac-b^2\right).$$
and we deduce the first part of Theorem \ref{geo}:
$$d\bar{A}= \sqrt{EG-F^2}dxdy =\sqrt{EG-F^2}e^{-2r} dA $$
$$=\Big( (w^2 + (2h +\Delta h) w + h^2+ h \Delta h + H (h) \Big)dA.$$

Next we calculate the coefficients of the second fundamental form:  since $N(u)=u,$ we have:
$$ l := \< \pa_x N(u), f_x\>=\< \pa_x u, f_x\>=e^{2r}(s+a),$$
$$ m:= \< \pa_x N(u), f_y\>=\< \pa_x u , f_y\>=e^{2r}b,$$
$$n :=\< \pa_y N(u), f_y\>=\< \pa_y u, f_y\>=e^{2r}(s+c).$$
Thus
$$lG+nE-2mF =e^{4r}\left((s+a)((s+c)^2+b^2)+(s+c)((s+a)^2+b^2)-2b^2(2s+a+c) \right)$$ 
$$=e^{4r}\left((s+a)(s+c)(2s+a+c)-b^2(2s+a+c) \right)$$
$$=e^{4r}(w^2+\al w + \be)(2w + \al),$$ 
and
$$ln-m^2=e^{4r}\left((s+c)(s+a)-b^2\right)=e^{4r}\left(w^2+\al w + \be \right).$$
Thus, at a point $f(u)$ where $d \bar{A}$ does not vanish, the mean curvature and the Gaussian curvature of $\pa B$
are given by
$$ 2H =
\frac{lG+nE-2mF}{EG-F^2}=\frac{2w + \al}{w^2 + \al w+ \be},$$
and
$$K=\frac{ln-m^2}{EG-F^2}=\frac{1}{w^2 +\al w + \be},$$
so that
$$H^2-K=\frac{(2w + \al)^2-4(w^2 +\al w + \be)}{4(w^2 +\al w + \be)^2}=\frac{\al^2 - 4\be}{4(w^2 +\al w + \be)^2}.$$
Hence the principal curvatures $k_1$ and $k_2$ of the immersion at the point $f(u)$ are 
$$ k_{1,2} = H \pm \sqrt{H^2-K}=\frac{2w + \al \pm \sqrt{\al^2-4 \be}}{2(w^2 + \al w+ \be)}.$$

\section{Volume and area of $\pa B$  (proof of Theorem~\ref{VA})}

The only tricky part of the proof is the following lemma:
\begin{lemm} \label{H} If $B$ has constant width, then
$$\int_{\S^2} H(h) dA= \frac{1}{2}\int_{\S^2} |\nabla h|^2 dA.$$
 \end{lemm}
 \noindent \textit{Proof.} Denoting the complex structure on $\S^2$ by $j$,  we have $j\pa_x = \pa_y, $  $j \pa_y =-\pa_x.$
 The proof is based on the following formula for the curvature tensor on the sphere:
 $$ \<R(X,Y)Z,W\> = \<X,Z\>\<Y,W\> - \<Y,Z\>\<X,W\>$$
 with $X=\pa_x, Y=\pa_y, Z=\nabla h, W = j\nabla h.$
 On the one hand,
 $$\int_{U} \big(\<X,Z\>\<Y,W\> - \<Y,Z\>\<X,W\> \big)dx dy $$
 $$=
  \int_{U} \left( \<\pa_x,\nabla h\> \<j\pa_x,j\nabla h\>-\<\pa_y,\nabla h\>\<\pa_y,j\nabla h\> \right)dx dy$$
  $$ = \int_{U} (h_x^2 + h_y^2)dx dy=\int_{\S^2} |\nabla h|^2 dA.$$
  On the other hand, using the fact that $j$ is parallel, i.e.\ $\nabla_X jY=j\nabla_X Y,$ we have
    $$\int_{U}  \<R(\pa_x,\pa_y )\nabla h,j \nabla h\> dx dy=
  \int_{U} (\<\nabla_{\pa_y} \nabla_{\pa_x} \nabla h,j \nabla h\>-\<\nabla_{\pa_x} \nabla_{\pa_y} \nabla h,j \nabla h\> )dx dy$$
  $$=\int_{U} (-\< \nabla_{\pa_x} \nabla h,\nabla_{\pa_y} j \nabla h\>+\< \nabla_{\pa_y} \nabla h, \nabla_{\pa_x}j \nabla h\> )dx dy$$
   $$=\int_{U} (-\< \nabla_{\pa_x} \nabla h,j\nabla_{\pa_y}  \nabla h\>+\< \nabla_{\pa_y} \nabla h,j \nabla_{\pa_x} \nabla h\> )dx dy$$
 $$=2 \int_{U}  \< \nabla_{\pa_y} \nabla h, j\nabla_{\pa_x} \nabla h\> )dx dy$$
 $$=2 \int_{U}  \< (b\pa_x+c\pa_y),(-b\pa_x+a\pa_y) \>dx dy$$
 $$=2 \int_{U}(ac-b^2)e^{2r} dx dy=2 \int_{\S^2} H(h)dA,$$
hence the proof of the lemma  is complete.

\bigskip

\bigskip

 In order to calculate the volume of $B$, we use the divergence theorem. Recalling that $u$ is
 the unit outward normal vector of the smooth parts of $\pa B,$ and that $f(u)=s(u)u + \nabla h (u),$
  we have
$${\cal V}(B)=\frac{1}{3}\int_{\S^2} \<f(u), u\> d \bar{A}=\frac{1}{3}\int_{\S^2} s(u) d \bar{A} $$
$$= \frac{1}{3}\int_{\S^2} (h+w)(w^2 + \al w + \be)  dA$$
$$=\frac{w^3}{3} \int_{\S^2} dA +\frac{w^2}{3} \int_{\S^2}(3h +\Delta h)dA+ 
\frac{w}{3} \int_{\S^2}(3h^2 + 2 h \Delta h+H(h))dA$$
$$+ \frac{1}{3} \int_{\S^2}(h^3 +h^2 \Delta h+ h H(h))dA.$$
Since $h$ has zero mean, the coefficient of $w^2$ vanishes. Moreover, the constant width condition, i.e.\ the oddness of
$h$,
implies that all the cubic expressions of $h$ and its second derivatives are odd and hence have
 zero mean. Thus the constant term vanishes.
 Finally, using the divergence theorem and Lemma \ref{H}, we obtain:
$${\cal V}(B)=\frac{w^3}{3} \int_{\S^2} dA + {w}\int_{\S^2}\left(h^2+\left(-\frac{2}{3}+\frac{1}{6}\right)|\nabla h|^2 \right)dA$$
$$=\frac{4 \pi w^3}{3}  - {w}\int_{\S^2}\left(\frac{1}{2}|\nabla h|^2-h^2 \right)dA$$
$$=\frac{4 \pi w^3}{3}  - {w} {\cal E}(h),$$
which is the required formula.

The computation of the area of $\pa B$ uses Lemma \ref{H} as well and is straightforward:
$${\cal A}(\pa B)=\int_{\pa B} d\bar{A}=\int_{\S^2}\left(w^2 +  (2h+\Delta h )w + (h^2+h\Delta h+ H(h))\right) dA $$
$$ =w^2 A(\S^2)  + \int_{\S^2} h^2 dA + \int_{\S^2}  |\nabla h|^2dA - \frac{1}{2}\int_{\S^2} |\nabla h |^2 dA$$
$$=4 \pi^2 - {\cal E}(h).$$


\bigskip

\noindent
Henri Anciaux \\
Universidade de S\~ao Paulo \\ 
  IME, Bloco A \\
  1010 Rua do Mat\~ao  \\
Cidade Universit\'aria   \\ 
 05508-090 S\~ao Paulo, Brazil   \\
henri.anciaux@gmail.com \\

\noindent
 Brendan Guilfoyle \\
Department of Mathematics and Computing \\
Institute of Technology, Tralee \\
Co. Kerry, Ireland \\
brendan.guilfoyle@ittralee.ie

\end{document}